%% file: admm_gpu.tex
\documentclass[a4paper,12pt]{article}

\input defs.tex

\title{%
  \bfseries\sffamily
    GPU Acceleration of ADMM for Large-Scale Quadratic Programming
}

\author{Michel Schubiger, Goran Banjac, and John Lygeros}

\newlength{\mywidth}
\newlength{\myheight}
\newcommand{\tabitem}{~--~}

\begin{document}

\maketitle

\begin{abstract}
  The \gls{ADMM} is a powerful operator splitting technique for solving structured convex optimization problems.
  Due to its relatively low per-iteration computational cost and ability to exploit sparsity in the problem data, it is particularly suitable for large-scale optimization.
  However, the method may still take prohibitively long to compute solutions to very large problem instances.
  Although \gls{ADMM} is known to be parallelizable, this feature is rarely exploited in real implementations.
  In this paper we exploit the parallel computing architecture of a \gls{GPU} to accelerate \gls{ADMM}.
  We build our solver on top of OSQP, a state-of-the-art implementation of \gls{ADMM} for quadratic programming.
  Our open-source CUDA C implementation has been tested on many large-scale problems and was shown to be up to two orders of magnitude faster than the CPU implementation.
\end{abstract}

\glsresetall

\section{Introduction}

Convex optimization has become a standard tool in many engineering fields including control \cite{Garcia:1989,Rawlings:2009}, signal processing \cite{Mattingley:2010}, statistics \cite{Huber:1964,Tibshirani:1996,Candes:2008}, finance \cite{Markowitz:1952,Cornuejols:2006,Boyd:2014,Boyd:2017}, and machine learning \cite{Cortes:1995}.
In some of these applications one seeks solutions to optimization problems whose dimensions can be very large.
For such problems, classical optimization algorithms, such as interior-point methods, may fail to provide a solution.

In the last decade operator splitting methods, such as the proximal gradient method and the \gls{ADMM}, have gained increasing attention in a wide range of application areas \cite{Boyd:2011,Parikh:2014,Banjac:2017}.
These methods scale well with the problem dimensions, can exploit sparsity in the problem data efficiently, and are often easily parallelizable.
Moreover, requirements on the solution accuracy are often moderate because of the noise in the data and arbitrariness of the objective.
This supports the use of operator splitting methods, which return solutions of a medium accuracy at a reasonable computational effort.

\Glspl{GPU} are hardware accelerators that offer an unmatched amount of parallel computational power for their relatively low price.
They provide far greater memory bandwidths than conventional CPU-based systems, which is especially beneficial in applications that process large amounts of data.
It is thus no surprise that the use of \glspl{GPU} has seen many applications in the area of machine learning, ranging from training deep neural networks \cite{Krizhevsky:2017,Ledig:2017,Cong:2019} to autonomous driving \cite{Liu:2017}.
Many software tools for machine learning, including PyTorch \cite{Steiner:2019}, TensorFlow \cite{TensorFlow:2016}, Theano \cite{Theano:2016}, and CNTK \cite{Seide:2016}, have native support for \gls{GPU} acceleration.
However, there has been a general perception that \glspl{GPU} are not well suited for the needs of numerical solvers for \glspl{LP} and \glspl{QP} \cite{gurobi-webinar}.

This paper explores the possibilities offered by the massive parallelism of \glspl{GPU} to accelerate solutions to large-scale \glspl{QP}.
We build our solver on top of the \gls{ADMM}-based OSQP solver \cite{Stellato:2020}.
The authors in \cite{O'Donoghue:2016} have demonstrated that \glspl{GPU} can be used to accelerate the solution to the linear system arising in their method.
We follow a similar approach to accelerate OSQP by replacing its direct linear system solver with an indirect (iterative) one, which we implement on the \gls{GPU}.
Moreover, we perform all vector and matrix operations on the \gls{GPU}, which further improves the performance of our implementation.
While the authors in \cite{Amos:2017,Charlton:2019} use \glspl{GPU} to solve \glspl{LP} and \glspl{QP} in batches, \ie they solve numerous different problems within one operation, our solver is designed for solving a single but large-scale problem at a time.

\subsection*{Outline}
We introduce the problem of interest in Section~\ref{sec:problem} and summarize the algorithm used by the OSQP solver in Section~\ref{sec:osqp}.
We then present in Section~\ref{sec:pcg} an alternative method for solving the linear system arising in OSQP.
We give a short summary of general GPU programming strategies in Section~\ref{sec:gpu}, followed by implementation details of the proposed \gls{GPU}-based solver in Section~\ref{sec:gpu_accel}.
Section~\ref{sec:numerics} demonstrates the performance of our solver on large-scale numerical examples.
Finally, Section~\ref{sec:conclusion} concludes the paper.

\subsection*{Notation}
Let $\Re$ denote the set of real numbers, $\Re^n$ the $n$-dimensional real space, $\Re^{m\times n}$ the set of real $m$-by-$n$ matrices, and $\symm_{++}^n$ ($\symm_+^n$) the set of real $n$-by-$n$ symmetric positive (semi)definite matrices.
We denote by $I$ and $\boldsymbol{1}$ the identity matrix and the vector of all ones (of appropriate dimensions), respectively.
For a vector $x\in\Re^n$, we denote its $i$-th element by $x_i$, the Euclidean norm by $\norm{x}_2 \eqdef \sqrt{x\tpose x}$, and the $\ell_\infty$ norm by $\norm{x}_\infty \eqdef \max_i |x_i|$.
For a matrix $K\in\symm_{++}^n$, we denote the $K$-norm of $x\in\Re^n$ by $\norm{x}_K \eqdef \sqrt{x\tpose Kx}$.
The gradient of a differentiable function $f\colon\Re^n\to\Re$ evaluated at $x\in\Re^n$ is denoted by $\nabla f(x)$.
For a nonempty, closed, and convex set $\mcf{C}\subseteq\Re^n$, we denote the Euclidean projection of $x\in\Re^n$ onto $\mcf{C}$ by $\Pi_\mcf{C}(x)\eqdef\argmin_{y\in\mcf{C}}\norm{x-y}_2$.
The Euclidean projections of $x\in\Re^n$ onto the nonnegative and nonpositive orthants are denoted by $x_+ \eqdef \max(x,0)$ and $x_- \eqdef \min(x,0)$, respectively.

\section{Problem Description}\label{sec:problem}

Consider the following \gls{QP}:
\begin{equation}
  \MinProblem[{eqn:qp}]{}{\half x\tpose P x + q\tpose x}{l \le Ax \le u,}
\end{equation}
where $x\in\Re^n$ is the optimization variable.
The objective function is defined by a positive semidefinite matrix $P \in \symm_+^n$ and a vector $q\in\Re^n$, and the constraints by a matrix $A\in\Re^{m\times n}$ and vectors $l$ and $u$ so that $l_i\in\Re\cup\{-\infty\}$, $u_i\in\Re\cup\{+\infty\}$, and $l_i \le u_i$ for all $i=1,\ldots,m$.
Linear equality constraints can be encoded in this way by setting $l_i=u_i$.

\subsection{Optimality and Infeasibility Conditions}

By introducing a variable $z\in\Re^m$, we can rewrite problem~\eqref{eqn:qp} in an equivalent form
\begin{equation}
  \MinProblem[{eqn:qpz}]{}{\half x\tpose P x + q\tpose x}{Ax = z, \quad l \le z \le u.}
\end{equation}
The optimality conditions for problem~\eqref{eqn:qpz} are given by \cite{Stellato:2020}
\begin{subequations}\label{eqn:opt_cond}
\begin{align}
  & Ax - z = 0 \label{eqn:prim_feas} \\
  & Px + q + A\tpose y = 0 \label{eqn:dual_feas}\\
  & l \le z \le u \label{eqn:box_constr}\\
  & y_-\tpose (z-l) = 0, \quad y_+\tpose (z-u) = 0,\label{eqn:compl_slack}
\end{align}
\end{subequations}
where $y\in\Re^m$ is a Lagrange multiplier associated with the constraint $Ax=z$.
If there exist $x\in\Re^n$, $z\in\Re^m$, and $y\in\Re^m$ that satisfy \eqref{eqn:opt_cond}, then we say that $(x,z)$ is a \emph{primal} and $y$ is a \emph{dual solution} to problem~\eqref{eqn:qpz}.

Problem~\eqref{eqn:qp} need not have a solution.
If there exists $\bar{y}\in\Re^m$ such that
\begin{equation}\label{eqn:prim_infeas}
  A\tpose \bar{y} = 0, \quad
  l\tpose \bar{y}_- + u\tpose \bar{y}_+ < 0,
\end{equation}
then problem~\eqref{eqn:qp} is infeasible and we say that $\bar{y}$ is a \emph{certificate of primal infeasibility}.
Similarly, if there exists $\bar{x}\in\Re^n$ such that
\begin{equation}\label{eqn:dual_infeas}
  P \bar{x} = 0, \quad
  q\tpose \bar{x} < 0, \quad
  (A \bar{x})_i \begin{cases}
		= 0   & l_i \in \Re, \, u_i \in \Re \\
		\ge 0 & u_i=+\infty \\
		\le 0 & l_i=-\infty
  \end{cases}
\end{equation}
for all $i=1,\ldots,m$, then the dual of problem~\eqref{eqn:qp} is infeasible and we say that $\bar{x}$ is a \emph{certificate of dual infeasibility}.
We refer the reader to \cite[\Prop~3.1]{Banjac:2019} for more details.

\section{OSQP Solver}\label{sec:osqp}

OSQP is an open-source numerical solver for convex \glspl{QP}.
It is based on ADMM and was shown to be competitive to and even faster than commercial \gls{QP} solvers \cite{Stellato:2020}.
An iteration of OSQP is shown in Algorithm~\ref{alg:osqp}.
The scalar $\alpha\in]0,2[$ is called the \emph{relaxation parameter}, and $\sigma>0$ and $R\in\symm_{++}^n$ are the \emph{penalty parameters}.
OSQP uses diagonal positive definite matrix $R$, which makes $R^{-1}$ easily computable.
In step~\ref{alg:update_z} of Algorithm~\ref{alg:osqp} it evaluates the Euclidean projection onto the box $[l,u]\eqdef\{z\in\Re^m \mid l \le z \le u\}$, which has a simple closed-form solution
\[
  \project{[l,u]}(z) = \min\left( \max\left( z, l \right), u \right),
\]
where $\min$ and $\max$ operators should be taken elementwise.

\begin{algorithm}[t]
  \caption{OSQP algorithm.}
  \label{alg:osqp}
  \begin{algorithmic}[1]
    \State \textbf{given} initial values $x^0$, $z^0$, $y^0$ and parameters $\sigma>0$, $R\in\symm_{++}^n$, $\alpha \in \left] 0,2 \right[$
    \State Set $k=0$
    \Repeat
    \State $(\tilde{x}^{k+1}, \tilde{z}^{k+1}) \gets \argmin\limits_{(\tilde{x},\tilde{z}): A\tilde{x}=\tilde{z}} \left\lbrace \half \tilde{x}\tpose P \tilde{x} + q\tpose \tilde{x} +
    \tfrac{\sigma}{2} \norm{\tilde{x} - x^k}_2^2 +
    \tfrac{1}{2} \norm{\tilde{z} - z^k + R\inv y^{k}}_R^2 \right\rbrace$ \label{alg:solve_lin_sys}
    \State $x^{k+1} \gets \alpha \tilde{x}^{k+1} + (1-\alpha)x^{k}$ \label{alg:update_x}
    \State $z^{k+1} \gets \project{[l,u]}\left(\alpha \tilde{z}^{k+1} + (1-\alpha)z^{k} + R\inv y^{k} \right)$ \label{alg:update_z}
    \State $y^{k+1} \gets y^{k} + R \left( \alpha \tilde{z}^{k+1} + (1-\alpha)z^{k} - z^{k+1} \right)$ \label{alg:update_y}
    \State $k \gets k+1$
    \Until{termination criterion is satisfied}
  \end{algorithmic}
\end{algorithm}

If problem~\eqref{eqn:qpz} is solvable, then the sequence $(x^k,z^k,y^k)$ generated by Algorithm~\ref{alg:osqp} converges to its primal-dual solution \cite{Banjac:2019,Stellato:2020}.
On the other hand, if the problem is primal or dual infeasible, then the iterates $(x^k,z^k,y^k)$ do not converge, but the sequence
\[
  (\delta x^k,\delta z^k,\delta y^k) \eqdef (x^k-x^{k-1},z^k-z^{k-1},y^k-y^{k-1})
\]
always converges and can be used to certify infeasibility of the problem.
In particular, if the problem is primal infeasible then $\delta y \eqdef \lim_{k\to\infty}\delta y^k$ will satisfy \eqref{eqn:prim_infeas}, whereas $\delta x \eqdef \lim_{k\to\infty}\delta x^k$ will satisfy \eqref{eqn:dual_infeas} if it is dual infeasible \cite[\Thm~5.1]{Banjac:2019}.

\subsection{Termination Criteria}\label{subsec:termination}
For the given iterates $(x^k,z^k,y^k)$, we define the primal and dual residuals as
\begin{align*}
  r_{\rm prim}^k &= Ax^k - z^k \\
  r_{\rm dual}^k &= Px^k + q + A\tpose y^k.
\end{align*}
The authors in \cite{Stellato:2020} show that the pair $(z^k,y^k)$ satisfies optimality conditions \eqref{eqn:box_constr}--\eqref{eqn:compl_slack} for all $k>0$ regardless of whether the problem is solvable or not.
If the problem is also solvable, then the residuals $r_{\rm prim}^k$ and $r_{\rm dual}^k$ will converge to zero \cite[\Prop~5.3]{Banjac:2019}.
A termination criterion for detecting optimality is thus implemented by checking that $r_{\rm prim}^k$ and $r_{\rm dual}^k$ are small enough, \ie
\begin{equation}\label{eqn:termination_condition}
  \norm{r_{\rm prim}^k}_\infty \le \eps_{\rm prim}, \quad
  \norm{r_{\rm dual}^k}_\infty \le \eps_{\rm dual},
\end{equation}
where $\eps_{\rm prim}>0$ and $\eps_{\rm dual}>0$ are some tolerance levels, which are often chosen relative to the scaling of the algorithm iterates \cite[\S 3.3]{Boyd:2011}.

Since $(\delta x^k, \delta y^k) \to (\delta x, \delta y)$, termination criteria for detecting primal and dual infeasibility are implemented by checking that $\delta y^k$ and $\delta x^k$ almost satisfy infeasibility conditions \eqref{eqn:prim_infeas} and \eqref{eqn:dual_infeas}, \ie
\begin{equation*}
	\norm{A\tpose \delta y^k}_\infty \le \eps_{\rm pinf}, \quad
	l\tpose (\delta y^k)_- + u\tpose (\delta y^k)_+ < \eps_{\rm pinf},
\end{equation*}
and
\begin{equation*}
  \norm{P \delta x^k}_\infty \le \eps_{\rm dinf}, \quad
  q\tpose \delta x^k < \eps_{\rm dinf}, \quad
  (A \delta x^k)_i \begin{cases}
    \in [-\eps_{\rm dinf},\eps_{\rm dinf}]	& l_i \in \Re, \, u_i \in \Re \\
    \ge -\eps_{\rm dinf}					& u_i=+\infty \\
    \le \eps_{\rm dinf}						& l_i=-\infty
  \end{cases}
\end{equation*}
for all $i=1,\ldots,m$, where $\eps_{\rm pinf}>0$ and $\eps_{\rm dinf}>0$ are given tolerance levels.

\subsection{Solving the KKT System}

Step~\ref{alg:solve_lin_sys} of Algorithm~\ref{alg:osqp} requires the solution to an equality-constrained \gls{QP}, which is equivalent to solving the following linear system:
\begin{equation}\label{eqn:kkt}
  \begin{bmatrix} P + \sigma I & A\tpose \\ A & -R\inv \end{bmatrix}
  \begin{bmatrix} \tilde{x}^{k+1} \\ \nu^{k+1} \end{bmatrix} =
  \begin{bmatrix} \sigma x^{k} - q \\ z^{k} - R\inv y^k \end{bmatrix},
\end{equation}
from which $\tilde{z}^{k+1}$ can be obtained as
\[
  \tilde{z}^{k+1} = z^k + R\inv (\nu^{k+1} - y^{k}).
\]
We refer to the matrix in \eqref{eqn:kkt} as the \emph{KKT matrix}.

OSQP uses a \emph{direct method} that computes the exact solution to \eqref{eqn:kkt} by first computing a factorization of the KKT matrix and then performing forward and backward substitutions.
The KKT matrix is symmetric quasi-definite for all $\sigma>0$ and $R\in\symm_{++}^n$, which ensures that it is nonsingular and has a well-defined $LDL\tpose$ factorization with diagonal $D$ \cite{George:2000}.
Since the KKT matrix does not depend on the iteration counter $k$, OSQP performs factorization at the beginning of the algorithm, and reuses the factors in subsequent iterations.

\subsection{Preconditioning}

A known weakness of \gls{ADMM} is its inability to deal effectively with ill-conditioned problems and its convergence can be very slow when data are badly scaled.
\emph{Preconditioning} is a common heuristic aiming to speed-up convergence of first-order methods.
OSQP uses a variant of \emph{Ruiz equilibration} \cite{Ruiz:2001,Knight:2014}, given in Algorithm~\ref{alg:ruiz}, which computes a cost scaling scalar $c>0$ and diagonal positive definite matrices $D$ and $E$ that effectively modify problem~\eqref{eqn:qp} into the following:
\[
  \MinProblem[{eqn:qp_mod}]{}{\half \bar{x}\tpose \bar{P} \bar{x} + \bar{q}\tpose \bar{x}}{\bar{l} \le \bar{A}\bar{x} \le \bar{u},}
\]
where the optimization variables are $\bar{x}=D\inv x$, $\bar{z}=E\inv z$ and $\bar{y}=cE\inv y$, and the problem data are
\[
  \bar{P}=cDPD, \quad
  \bar{q}=cDq,  \quad
  \bar{A}=EAD,  \quad
  \bar{l}=El,   \quad
  \bar{u}=Eu.
\]

\begin{algorithm}[t]
  \caption{Modified Ruiz equilibration.}
  \label{alg:ruiz}
  \begin{algorithmic}[1]
    \State \textbf{initialize} $c=1$, $D=I$, $E=I$, $\delta = 0$, $\bar{P}=P$, $\bar{q}=q$, $\bar{A}=A$
    \While{$\norm{1 - \delta}_\infty > \eps_{\rm equil}$}
      \State $M \gets \begin{bsmallmatrix} \bar{P} & \bar{A}\tpose \\ \bar{A} & 0 \end{bsmallmatrix}$
      \For{$i = 1,\dots,n+m$}
        \State $\delta_i \gets 1 / \sqrt{\norm{M_i}_{\infty}}$ \label{alg:column_norm}
      \EndFor
      \State $\begin{bsmallmatrix} D & \\ & E \end{bsmallmatrix} \gets \diag(\delta) \begin{bsmallmatrix} D & \\ & E \end{bsmallmatrix}$
      \State $\bar{P} \gets D P D, \quad \bar{q} \gets D q, \quad \bar{A} \gets E A D$ \label{alg:matrix_scaling}
      \State $\gamma \gets 1/\max\lbrace\mean(\norm{\bar{P}_i}_{\infty}), \norm{\bar{q}}_{\infty}\rbrace$
      \State $\bar{P} \gets \gamma \bar{P}, \quad \bar{q} \gets \gamma \bar{q}, \quad c \gets \gamma c$
    \EndWhile
    \Return{$D$, $E$, $c$}
  \end{algorithmic}
\end{algorithm}

\subsection{Parameter Selection}

OSQP sets $\alpha=1.6$ and $\sigma=10^{-6}$ by default and the choice of these parameters does not seem to be critical for the \gls{ADMM} convergence rate.
However, the choice of $R=\diag(\rho_1,\ldots,\rho_m)$ is a key determinant of the number of iterations required to satisfy a termination criterion.
OSQP sets a higher value of $\rho_i$ that is associated with an equality constraint, \ie
\begin{equation}\label{eqn:rho_vec}
  \rho_i=\begin{cases} \bar{\rho} & l_i \ne u_i \\ 10^3 \bar{\rho} & l_i = u_i \end{cases},
\end{equation}
where $\bar{\rho}>0$.
Having a fixed value of $\bar{\rho}$ does not provide satisfactory performance of the algorithm across different problems.
To compensate for this sensitivity, OSQP adopts an adaptive scheme, which updates $\bar{\rho}$ during the iterations based on the ratio between norms of the primal and dual residuals \cite{Wohlberg:2017}.

The proposed parameter update scheme makes the algorithm much more robust, but also introduces additional computational burden since updating $R$ changes the KKT matrix in \eqref{eqn:kkt}, which then needs to be refactored.
Updating $\bar{\rho}$ is thus performed only a few times during the runtime of the algorithm.

\section{Preconditioned Conjugate Gradient Method}\label{sec:pcg}

An alternative way to solve the equality-constrained \gls{QP} in step~\ref{alg:solve_lin_sys} of Algorithm~\ref{alg:osqp} is by using an \emph{indirect method}.
As observed in \cite{Stellato:2020}, eliminating $\nu^{k+1}$ from \eqref{eqn:kkt} results in the \emph{reduced KKT system}
\begin{equation}\label{eqn:kkt_red}
  (P + \sigma I + A\tpose R A) \tilde{x}^{k+1} = \sigma x^k - q + A\tpose (R z^k - y^k),
\end{equation}
from which $\tilde{z}^{k+1}$ can be obtained as $\tilde{z}^{k+1}=A\tilde{x}^{k+1}$.
Note that the reduced KKT matrix is always positive definite, which allows us to use the \gls{CG} method for solving~\eqref{eqn:kkt_red}.

\subsection{Conjugate Gradient Method}

The \gls{CG} method is an iterative method for solving linear systems of the form
\begin{equation}\label{eqn:pd_lin_sys}
  Kx=b,
\end{equation}
where $K\in\symm^n_{++}$ is a symmetric positive definite matrix.
The method computes the solution to the linear system in at most $n$ iterations \cite[\Thm~5.1]{Nocedal:2006}.
However, when solving large-scale linear systems, one aims to terminate the method after $d \ll n$ iterations, which yields an approximate solution to~\eqref{eqn:pd_lin_sys}.

Solving \eqref{eqn:pd_lin_sys} is equivalent to solving the following unconstrained optimization problem:
\[
  \text{minimize} \quad f(x) \eqdef \half x\tpose K x - b\tpose x,
\]
since its minimizer can be characterized as
\[
  0 = \nabla f(x) = Kx-b.
\]

\subsubsection{Conjugate directions}
A set of nonzero vectors $\{p^0,\ldots,p^{n-1}\}$ is said to be \emph{conjugate} with respect to $K$ if
\[
  (p^i)\tpose K p^j = 0, \quad \forall i \ne j.
\]
Successive minimization of $f$ along the conjugate directions $p^k$, \ie evaluating
\begin{subequations}\label{eqn:succ_min}
\begin{align}
  \alpha^k &= \argmin_\alpha f(x^k + \alpha p^k) \label{eqn:succ_min:alpha} \\
  x^{k+1}	 &= x^k + \alpha^k p^k,
\end{align}
\end{subequations}
produces $x^{k+1}$ that minimizes $f$ over $(\{x^0\} + S^k)$, where $S^k$ is the expanding subspace spanned by the previous conjugate directions $\{p^0,\ldots,p^k\}$ \cite[\Thm~5.2]{Nocedal:2006}.
The minimization in~\eqref{eqn:succ_min:alpha} has the following closed-form solution:
\[
  \alpha^k = -\frac{(r^k)\tpose p^k}{(p^k)\tpose K p^k},
\]
where $r^k \eqdef \nabla f(x^k) = Kx^k-b$ is the residual at step~$k$.

\subsubsection{Conjugate gradient}
There are various choices for the conjugate direction set $\{p^0,\ldots,p^{n-1}\}$.
For instance, the eigenvectors of $K$ form a set of conjugate directions with respect to $K$, but are impractical to compute for large matrices.
The cornerstone of the \gls{CG} method is its ability to generate a set of conjugate directions efficiently.
It computes a new direction $p^k$ using only the previous direction $p^{k-1}$, which imposes low computational and memory requirements.
In particular, a new direction $p^k$ is computed as a linear combination of the negative gradient $-r^k$ and the previous direction $p^{k-1}$,
\begin{equation}\label{eqn:new_conj_dir}
  p^k = -r^k + \beta^k p^{k-1}.
\end{equation}
The scalar $\beta^k$ is determined from the conjugacy requirement $(p^k)\tpose Kp^{k-1} = 0$, leading to
\[
  \beta^k = \frac{(r^k)\tpose Kp^{k-1}}{(p^{k-1})\tpose Kp^{k-1}}.
\]
The first conjugate direction is set to the negative gradient, \ie $p^0=-r^0$.
Combining the successive minimization~\eqref{eqn:succ_min} and the computation of conjugate directions~\eqref{eqn:new_conj_dir} yield the \gls{CG} method.

\subsection{Preconditioning}
Since the \gls{CG} method is a first-order optimization method, it is sensitive to the problem scaling.
To improve convergence of the method, we can precondition the linear system by using a coordinate transformation
\[
  \bar{x} = Cx,
\]
where $C\in\Re^{n\times n}$ is a nonsingular matrix.
Applying the \gls{CG} method to the transformed linear system yields the \gls{PCG} method, which is shown in Algorithm~\ref{alg:pcg} \cite[\Alg~5.3]{Nocedal:2006}.
It turns out that $C$ need not be formed explicitly, but rather acts through $M=C\tpose C$.

\begin{algorithm}[t]
  \caption{\Gls{PCG} method.}
  \label{alg:pcg}
  \begin{algorithmic}[1]
    \Require{initial value $x^0$, preconditioner $M$}
    \State \textbf{initialize} $r^0 = K x^0 - b$, $y^0 = M^{-1}r^0$, $p^0 = -y^0$
    \State Set $k=0$
    \While{$\norm{r^k} > \eps \norm{b}$} \label{alg:pcg_termination}
      \State $\alpha^k \gets -\frac{(r^k)\tpose y^k}{(p^k)\tpose K p^k}$
      \State $x^{k+1} \gets x^k + \alpha^k p^k$
      \State $r^{k+1} \gets r^k + \alpha^k K p^{k}$
      \State $y^{k+1} \gets M^{-1}r^{k+1}$
      \State $\beta^{k+1} \gets \frac{(r^{k+1})\tpose y^{k+1}}{(r^k)\tpose y^{k}}$
      \State $p^{k+1} \gets -y^{k+1} + \beta^{k+1} p^k$
      \State $k \gets k+1$
    \EndWhile
  \end{algorithmic}
\end{algorithm}

In general, a good preconditioner should satisfy $M \approx K$ and at the same time make the linear system $My=r$ easy to solve \cite[\S 11.5]{Golub:2013}.
One of the simplest choices is the \emph{diagonal} or \emph{Jacobi preconditioner}, which contains the diagonal elements of $K$, making $M\inv$ easily computable.

More advanced choices include the incomplete Cholesky, the incomplete $LU$, and polynomial preconditioners.
The incomplete preconditioners produce an approximate decomposition of $K$ with a high sparsity, so that solving $My=r$ is computationally cheap.
The family of polynomial preconditioners include the Chebyshev and least-squares polynomial preconditioners, both of which require a bound on the spectrum of $K$ \cite{Ashby:1992}.

\section{GPU Architecture and Programming Strategies}\label{sec:gpu}

\glspl{GPU} have been used for general-purpose computing for more than two decades \cite{Brodtkorb:2013}.
They come in many different variations and architectures, but we will restrict our discussion to the latest NVIDIA Turing-based architecture.
Most of the concepts also apply to older NVIDIA \glspl{GPU}; for further details, we refer the reader to~\cite{cuda}.

A \gls{GPU} consists of an array of several \glspl{SM}, each of which contains multiple integer and floating-point arithmetic units, local caches, shared memory, and several schedulers.
The on-board RAM of the GPU is called \emph{global memory}, in contrast to the shared memory that is local to each \gls{SM}.
While an \gls{SM} of the Turing generation has $96~\rm{KB}$ of shared memory, the global memory is much larger and is typically in order of $\rm{GB}$.
Compared to the shared memory, it has a much higher latency and a lower bandwidth, but is still much faster than system RAM; bandwidths of $500~\rm{GB/s}$ are not uncommon for the \gls{GPU} global memory, whereas system RAM is limited to $40$-$50~\rm{GB/s}$.

The main challenge in using \glspl{GPU} is to leverage the increasing number of processing cores and develop applications that scale their parallelism.
A solution designed by NVIDIA to overcome this challenge is called CUDA, a general-purpose parallel computing platform and programming model.

\subsection{CUDA Architecture}
CUDA is an extension of the C programming language created by NVIDIA.
Its main idea is to have a large number of threads that solve a problem cooperatively.
This section explains how threads are organized into cooperative groups and how CUDA achieves scalability.

\subsubsection{Kernels}
Kernels are special C functions that are executed on a \gls{GPU} and are defined by the \texttt{\_\_global\_\_} keyword.
In contrast to regular C functions, kernels get executed $N$ times in parallel by $N$ different threads, where each thread executes the same code, but on different data.
This concept is known as \gls{SIMD}.
The number of threads is specified when calling a kernel, which is referred to as a kernel launch.

\subsubsection{Thread hierarchy}
While kernels specify the code that is executed by each thread, the thread hierarchy dictates how the individual threads are organized.
CUDA has a two-level hierarchy to organize the threads, a \emph{grid-level} and a \emph{block-level}.
A grid contains multiple blocks and a block contains multiple threads.
A kernel launch specifies the grid size (number of blocks) and the block size (number of threads per block).

The threads within one block can cooperate to solve a subproblem.
The problem needs to be partitioned into independent subproblems by the programmer so that a grid of thread blocks can solve it in parallel.
Each block is scheduled on one of the available \glspl{SM}, which can happen concurrently or sequentially, depending on the number of blocks and available hardware.
If there are enough resources available, then several blocks can be scheduled on a single \gls{SM}.

The threads within a single block have a unique thread index that is accessible through the built-in variable \texttt{threadIdx}, which is defined as a $3$-dimensional vector.
This allows threads to be indexed in one-, two-, or three-dimensional blocks, which allows for a natural indexing in the problem domain.
Similarly, blocks within a grid have a unique block index that is accessible through the variable \texttt{blockIdx}.
It is also defined as a $3$-dimensional vector and allows for one-, two-, or three-dimensional indexing.

\subsubsection{Accelerating numerical methods}
Numerical methods make extensive use of floating-point operations, but their performance is not solely determined by the system's floating-point performance.
Although \glspl{GPU} offer magnitudes larger floating-point power than CPUs, it is the memory bandwidth that limits the performance of many numerical operations \cite{Vuduc:2010}.
Fortunately, \glspl{GPU} also offer an order of magnitude larger memory bandwidth, but utilizing their full potential is not an easy task since the parallel nature of the \gls{GPU} requires different programming strategies.

Listing~\ref{lst:axpy_cpu} implements the simple \gls{BLAS} routine \texttt{axpy} ($y=y+ax$) on the CPU.
The code uses a simple \texttt{for} loop to iterate through the elements of $x$ and $y$.
A \gls{GPU} implementation of \texttt{axpy} is shown in Listing~\ref{lst:axpy_gpu}.
The code looks very similar to the CPU version, but has two important differences.
First, the \texttt{for} loop is replaced by a simple \texttt{if} condition.
Instead of one thread iterating through a loop element-by-element, there is a thread for each element to be processed, which is a common pattern in \gls{GPU} computing.
Second, a thread ID is used to determine the data element which each thread is operating on.
A global thread ID is calculated from a local thread index and a block index.
The \texttt{if} condition disables threads with a thread ID larger than the number of elements.
This is necessary since threads are launched in blocks and the total number of threads usually does not match the number of elements, while the cost of few idle threads is negligible.

\begin{lstlisting}[
  float,
  floatplacement=t,
  language=C,
  caption={A CPU implementation of the \texttt{axpy} function.},
  label={lst:axpy_cpu}
]
  void axpy_cpu(float* y, float* x, float a, int n) {
    for (int idx = 0; idx < n; idx++) {
      y[idx] += a * x[idx];
    }
  }
\end{lstlisting}

\begin{lstlisting}[
  float,
  floatplacement=t,
  language=C,
  caption={A GPU implementation of the \texttt{axpy} function.},
  label={lst:axpy_gpu}
]
  __global__ void axpy_gpu(float* y, float* x, float a, int n) {
    int idx = threadIdx.x + blockDim.x * blockIdx.x;
    if (idx < n) {
      y[idx] += a * x[idx];
    }
  }
\end{lstlisting}

Figure~\ref{fig:axpy} compares the achieved memory throughput and the floating-point performance of the CPU and \gls{GPU} implementations of the \texttt{axpy} operation; see Section~\ref{sec:numerics} for the hardware specifications.
The plots are obtained by performing the \texttt{axpy} operation for various sizes of vectors and measuring the time required to run these operations.
Knowing how much data is moved and how many floating-point operations are required for the \texttt{axpy} operation, we can compute the average memory throughput and the average floating-point performance, both of which depend linearly on the size of vectors; this is why the shapes of the two plots look the same.
For large vector sizes the simple \gls{GPU} implementation is approximately $15$ times faster.
Note, however, that small problems cannot be accelerated well with \glspl{GPU}, as there is not enough work to keep the \gls{GPU} busy and a kernel launch and data transfer come with a constant overhead that cannot be amortized.
The \gls{GPU} reaches the maximum memory throughput of $522~\rm{GB/s}$, which is $85\%$ of its theoretical peak of $616~\rm{GB/s}$, whereas the peak value of the floating-point performance is around $87~\rm{GFLOPS}$, which is less than $1\%$ of its theoretical peak of $13.45~\rm{TFLOPS}$; the peak numbers are specifications of the NVIDIA GeForce RTX 2080 Ti \gls{GPU}, which was used in our numerical tests.
This shows that the performance of \texttt{axpy} is clearly limited by the memory bandwidth.

\begin{figure}[t]
  \centering
  \footnotesize
  \begin{tikzpicture}
    \begin{axis}[%
      hide axis,
      xmin=10,	
      xmax=50,	
      ymin=0,		
      ymax=0.4,	
      legend cell align=left,
      legend columns=3,
      legend style={%
        nodes={scale=1},
        fill=white,
        fill opacity=1,
        draw opacity=1,
        text opacity=1,
        /tikz/every even column/.append style={column sep=3em}
      },
      every axis plot/.append style={very thick}
      ]
      \addlegendimage{blue, mark=*, mark options={scale=1, fill=blue}}
      \addlegendentry{CPU (\texttt{-O0})};
      \addlegendimage{magenta, mark=pentagon*, mark options={scale=1, fill=magenta}}
      \addlegendentry{CPU (\texttt{-O3})};
      \addlegendimage{red, mark=square*, mark options={scale=1, fill=red}}
      \addlegendentry{GPU};
    \end{axis}
  \end{tikzpicture}
  \\[.5em]
  \begin{tabular}{rl}
    \begin{tikzpicture}
      \begin{loglogaxis}[%
        width=\mywidth,
        height=\myheight,
        xlabel={Problem size $n$},
        grid=both,
        grid style={line width=.1pt, draw=gray!10},
        major grid style={line width=.2pt,draw=gray!50},
        ]
        \addplot [blue, mark=*, mark options={scale=1, fill=blue}] table[x=N, y=CPU_Bandwidth, col sep=comma] {data/axpy.csv};
        \addplot [magenta, mark=pentagon*, mark options={scale=1, fill=magenta}] table[x=N, y=CPU_O3_Bandwidth, col sep=comma] {data/axpy.csv};
        \addplot [red, mark=square*, mark options={scale=1, fill=red}] table[x=N, y=GPU_Bandwidth, col sep=comma] {data/axpy.csv};
      \end{loglogaxis}
    \end{tikzpicture} 
    &
    \begin{tikzpicture}
      \begin{loglogaxis}[%
        width=\mywidth,
        height=\myheight,
        xlabel={Problem size $n$},
        grid=both,
        grid style={line width=.1pt, draw=gray!10},
        major grid style={line width=.2pt,draw=gray!50},
        ylabel near ticks,
        yticklabel pos=right,
        ]
        \addplot [blue, mark=*, mark options={scale=1, fill=blue}] table[x=N, y=CPU_GFLOPS, col sep=comma] {data/axpy.csv};
        \addplot [magenta, mark=pentagon*, mark options={scale=1, fill=magenta}] table[x=N, y=CPU_O3_GFLOPS, col sep=comma] {data/axpy.csv};
        \addplot [red, mark=square*, mark options={scale=1, fill=red}] table[x=N, y=GPU_GFLOPS, col sep=comma] {data/axpy.csv};
      \end{loglogaxis}
    \end{tikzpicture}
  \end{tabular}
  \caption{%
    Numerical performance of the \texttt{axpy} routine run on CPU and \gls{GPU}.
    Left: The average memory throughput (in $\rm{GB/s}$).
    Right: The average floating-point performance (in $\rm{GFLOPS}$).
  }
  \label{fig:axpy}
\end{figure}

\subsubsection{Segmented reduction}
A reduction is an operation that takes a vector $x\in\Re^n$ and an associative binary operator $\oplus$, and returns a scalar $y\in\Re$ \cite{Lammel:2008},
\[
	y = x_1 \oplus x_2 \oplus \ldots \oplus x_n.
\]
This abstract formulation allows us to formulate many operations as a reduction, among others the sum of elements, the maximum value of elements, the $\ell_1$ norm, the $\ell_\infty$ norm etc.
The only difference between reduction and segmented reduction is that the latter reduces individual segments of $x$ and outputs a vector that computes reduction over the segments.
There exist very efficient parallel implementations for both reduction and segmented reduction \cite{moderngpu,thrust}, and thus any problem that can be reformulated as one of them can be easily accelerated by a \gls{GPU}.

\subsection{CUDA Libraries}
There exist multiple libraries shipped with the CUDA Toolkit that implement various functions on the \gls{GPU} \cite{cuda-toolkit}.
We summarize in the sequel the NVIDIA libraries used in this work.

\begin{itemize}
  \item \emph{Thrust} is a CUDA \Cpp library based on the \Cpp \gls{STL}.
  It provides a high-level interface for high-performance parallel applications and all essential data parallel primitives, such as scan, sort, and reduce.
  
  \item \emph{cuBLAS} is a CUDA implementation of \gls{BLAS}, which enables easy \gls{GPU} acceleration of code that uses \gls{BLAS} functions.
  We use only \mbox{level-1} cuBLAS API functions that implement the inner product, \texttt{axpy} operation, scalar-vector multiplication, and computation of norms.
  
  \item \emph{cuSPARSE} is a CUDA library that contains a set of linear algebra subroutines for handling sparse matrices.
  It requires the matrices to be in one of the sparse matrix formats described in the next section.
\end{itemize}

\subsection{Sparse Matrix Formats}

\subsubsection{COO}
The \gls{COO} format is one of the simplest sparse matrix formats.
It is mainly used as an intermediate format to perform matrix operations, such as transpose, concatenation, or the extension of an upper triangular to a full symmetric matrix.
It holds the number of rows \texttt{m}, the number of columns \texttt{n}, the number of nonzero elements \texttt{nnz}, and three arrays of dimension \texttt{nnz}: \texttt{Value}, \texttt{RowIndex}, and \texttt{ColumnIndex}.
The cuSPARSE API assumes that the indices are sorted by rows first and then by columns within each row, which makes the representation unique.

The $4\times 5$ matrix given below:
\begin{equation}\label{eqn:A_mat_format}
  A = \begin{bmatrix} 1 & 0 & 0 & 0 & 4 \\ 0 & 5 & 1 & 0 & 0 \\ 0 & 2 & 0 & 0 & 1 \\ 7 & 0 & 1 & 0 & 0 \end{bmatrix}
\end{equation}
has the following \gls{COO} representation:
\begin{align*}
  \texttt{Value}       &= \begin{bmatrix} 1 & 4 & 5 & 1 & 2 & 1 & 7 & 1 \end{bmatrix} \\
  \texttt{RowIndex}    &= \begin{bmatrix} 0 & 0 & 1 & 1 & 2 & 2 & 3 & 3 \end{bmatrix} \\
  \texttt{ColumnIndex} &= \begin{bmatrix} 0 & 4 & 1 & 2 & 1 & 4 & 0 & 2 \end{bmatrix}.
\end{align*}
Note that we use the zero-based indexing in the example above and throughout the paper.

\subsubsection{CSR}
The \gls{CSR} format differs from the \gls{COO} format only in the \texttt{RowIndex} array, which is compressed in the \gls{CSR} format.
The compression can be understood as a two-step process.
First, we determine from \texttt{RowIndex} the number of nonzero elements in each row, which results in an array of length \texttt{m}.
Then, we calculate the cumulative sum of this array and insert a zero at the beginning, which results in an array of length \texttt{m+1}.
The obtained array is denoted by \texttt{RowPointer} since it points to the beginning of a row in \texttt{Value} and \texttt{ColumnIndex} arrays.

The \texttt{RowPointer} array has the property that the difference between its two consecutive elements,
\[
  \texttt{RowPointer[k+1] - RowPointer[k]},
\]
is equal to the number of nonzero elements in row \texttt{k}.
Noting that \texttt{RowPointer[0] = 0} and applying the property above recursively, it follows that
\[
  \texttt{RowPointer[m] = nnz}.
\]
Matrix $A$ given in~\eqref{eqn:A_mat_format} has the following \gls{CSR} representation:
\begin{align*}
  \texttt{Value}       &= \begin{bmatrix} 1 & 4 & 5 & 1 & 2 & 1 & 7 & 1 \end{bmatrix} \\
  \texttt{RowPointer}	 &= \begin{bmatrix} 0 & 2 & 4 & 6 & 8 \end{bmatrix} \\
  \texttt{ColumnIndex} &= \begin{bmatrix} 0 & 4 & 1 & 2 & 1 & 4 & 0 & 2 \end{bmatrix}.
\end{align*}
The \gls{CSR} format is used for \gls{SpMV} in cuSPARSE.

\subsubsection{CSC}
The \gls{CSC} format differs from the \gls{CSR} format in two ways: the values are stored in the column-major format and the column indices are compressed.
The compressed array has dimension \texttt{n+1} and is denoted by \texttt{ColumnPointer}.
Matrix $A$ given in~\eqref{eqn:A_mat_format} has the following \gls{CSC} representation:
\begin{align*}
  \texttt{Value}         &= \begin{bmatrix} 1 & 7 & 5 & 2 & 1 & 1 & 4 & 1 \end{bmatrix} \\
  \texttt{RowIndex}      &= \begin{bmatrix} 0 & 3 & 1 & 2 & 1 & 3 & 0 & 2 \end{bmatrix} \\
  \texttt{ColumnPointer} &= \begin{bmatrix} 0 & 2 & 4 & 6 & 6 & 8 \end{bmatrix}.
\end{align*}
The \gls{CSC} format is not used directly for computations in cuSPARSE.
However, we can interpret the \gls{CSC} representation of a matrix $A$ as the \gls{CSR} representation of $A\tpose$ using the following mapping:
\begin{align*}
  \texttt{m}_\texttt{CSC}             &\to \texttt{n}_\texttt{CSR} \\
  \texttt{n}_\texttt{CSC}             &\to \texttt{m}_\texttt{CSR} \\
  \texttt{ColumnPointer}_\texttt{CSC} &\to \texttt{RowPointer}_\texttt{CSR} \\
  \texttt{RowPointer}_\texttt{CSC}    &\to \texttt{ColumnPointer}_\texttt{CSR} \\
  \texttt{Value}_\texttt{CSC}         &\to \texttt{Value}_\texttt{CSR}.
\end{align*}

\section{GPU Acceleration of OSQP}\label{sec:gpu_accel}

Profile-driven software development is based on identifying major computational bottlenecks in the code, as performance will increase the most when removing these \cite{Brodtkorb:2013}.
This section identifies and analyzes computational bottlenecks of OSQP when solving large-scale \glspl{QP}, and shows how we can remove them by making use of \gls{GPU}'s parallelism.

\subsection{OSQP Computational Bottlenecks}\label{subsec:bottlenecks}
Given a \gls{QP} in the form \eqref{eqn:qp}, we denote the total number of nonzero elements in matrices $P$ and $A$ by $N\eqdef\nnz(P)+\nnz(A)$.
Profiling the OSQP code reveals that for large-scale problem instances all operations whose execution time scales with $N$ represent a potential bottleneck since $N$ is typically much larger than the number of \gls{QP} variables $n$ and constraints $m$.

As shown in Section~\ref{subsec:termination}, evaluating termination criteria requires several sparse matrix-vector multiplications.
Performing these computations in each ADMM iteration can slow down the solver considerably.
Hence, OSQP evaluates these criteria every $25$ iterations by default so that the overall computational burden is reduced.
This means that the algorithm can terminate only when the iteration counter $k$ is a multiple of $25$.
We discuss in Section~\ref{subsec:matrix-representation} how to represent the problem matrices in the GPU memory so that sparse matrix-vector multiplications can be performed efficiently on the GPU.

The main computational bottleneck is using a direct linear system solver to tackle the KKT system \eqref{eqn:kkt}.
When $N$ is very large, the computational cost of factoring the KKT matrix becomes prohibitively large.
This issue also limits the number of parameter updates, which can improve convergence rate of the algorithm, but require the KKT matrix to be refactored.
Furthermore, in each ADMM iteration we need to evaluate forward and backward substitutions, which cannot be fully parallelized.
Section~\ref{subsec:reduced-kkt} describes an efficient GPU implementation of the \gls{PCG} method that avoids matrix factorizations.

Profiling reveals that the matrix equilibration procedure described in Algorithm~\ref{alg:ruiz} is also demanding for large-scale problems, where the main bottlenecks are computing the column-norms in step~\ref{alg:column_norm} and matrix scaling in step~\ref{alg:matrix_scaling}.
The matrix scaling requires pre- and post-multiplying $P$ and $A$ by diagonal matrices, which is equivalent to scaling rows and columns of $P$ and $A$.
We discuss in Section~\ref{subsec:matrix_equilibration} how to parallelize these operations on the GPU.

\subsection{Representation of Matrices}\label{subsec:matrix-representation}
OSQP represents matrices $P$ and $A$ in the \gls{CSC} format.
Moreover, since $P$ is symmetric, only the upper triangular part of $P$ is actually stored in memory.
The preferred way of storing matrices in the \gls{GPU} memory is using the \gls{CSR} format since it has a superior \gls{SpMV} performance on the \gls{GPU}.
However, when using $A$ in the \gls{CSR} format, computing $A\tpose y$ is around $10$ times slower than computing $Ax$ \cite{cuda-toolkit}.
This inefficiency can be avoided by storing both $A$ and $A\tpose$ in the \gls{CSR} format, though this doubles the memory requirements.

Similarly, storing only the upper triangular part of $P$ is memory-efficient, but computing $Px$ in that case is much slower than when the full $P$ is stored \cite{cuda-toolkit}.
Therefore, we store the full $P$ in the CSR format since it improves the \gls{SpMV} performance.

We also store vectors $q$, $l$, $u$, as well as the ADMM iterates in the GPU memory, and perform all matrix and vector operations on the GPU.
This reduces considerably the size of memory transferred between the system and the GPU memory.

\subsection{Reduced KKT System}\label{subsec:reduced-kkt}
As discussed in Section~\ref{sec:pcg}, we can avoid factoring the KKT matrix by solving the reduced KKT system \eqref{eqn:kkt_red} with the \gls{PCG} method, which only evaluates matrix-vector multiplications and can be easily parallelized.
Although OSQP uses the parameter matrix of the form $R=\diag(\rho_1,\ldots,\rho_m)$, where $\rho_i$ is set as in \eqref{eqn:rho_vec}, numerical tests show that this choice of $R$ makes the \gls{PCG} method converge slowly.
This can be understood by looking at the effect of $R$ on the reduced KKT matrix.
Since $R$ appears in the term $A\tpose R A$, setting it as in \eqref{eqn:rho_vec} has the effect of scaling the rows of $A$ by different values, which effectively increases the condition number of the matrix.

The convergence rate of the \gls{PCG} method can be improved by using $R=\bar{\rho}I$ instead.
This choice will in general result in more iterations of Algorithm~\ref{alg:osqp}, but will reduce the number of iterations of Algorithm~\ref{alg:pcg} considerably.
The linear system~\eqref{eqn:kkt_red} now reduces to
\begin{equation}\label{eqn:kkt_red_sc}
  (P + \sigma I + \bar{\rho} A\tpose A) \tilde{x}^{k+1} = \sigma x^k - q + A\tpose (\bar{\rho} z^k - y^k).
\end{equation}
Note that the coefficient matrix above need not be formed explicitly.
Instead, the matrix-vector product
\[
  r \gets (P + \sigma I + \bar{\rho} A\tpose A)x
\]
can be evaluated as
\begin{align*}
  z &\gets \bar{\rho} Ax \\
  r &\gets Px + \sigma x + A\tpose z.
\end{align*}

\subsubsection{Preconditioner}
We use the Jacobi preconditioner, for which solving $My=r$ amounts to a simple diagonal matrix-vector product.
The diagonal of the Jacobi preconditioner for~\eqref{eqn:kkt_red_sc} can be computed as
\[
  \diag(M) = \diag(P) + \sigma \boldsymbol{1} + \bar{\rho} \diag(A\tpose A).
\]
Note that we need not compute the full product $A\tpose A$, but only its diagonal elements,
\[
  (A\tpose A)_{ii} = \norm{A_i}_2^2,
\]
where $A_i$ denotes the $i$-th column of $A$.

\subsubsection{Parameter update}
Once $\diag(P)$ and $\diag(A\tpose A)$ are available, computing $M$ becomes extremely easy.
This makes the parameter update computationally cheap since we only need to update the preconditioner $M$.
This allows us to update $\bar{\rho}$ more often than is done in OSQP.
Our numerical tests perform well when $\bar{\rho}$ is updated every $10$ iterations.

\subsubsection{Termination criteria and warm starting}
The solution to~\eqref{eqn:kkt_red_sc} need not be carried out exactly for Algorithm~\ref{alg:osqp} to converge \cite[\S 3.4.4]{Boyd:2011}.
This fact can be used to motivate an early termination of the \gls{PCG} method, meaning that we solve \eqref{eqn:kkt_red_sc} only approximately at first, and then more accurately as the iterations progress.
This can be achieved by performing a relatively small number of \gls{PCG} iterations to obtain an approximate solution, and using warm-starting by initializing $x^0$ in Algorithm~\ref{alg:pcg} to the solution $\tilde{x}^k$ computed in the previous \gls{ADMM} iteration.

Finding a good termination criterion for Algorithm~\ref{alg:pcg} is essential for reducing the total runtime of \gls{ADMM}.
If the \gls{PCG} method returns solutions with low accuracy, then \gls{ADMM} may converge slower, or even diverge.
On the other hand, if the \gls{PCG} method solves the subproblems with unnecessarily high accuracy, this may increase the total runtime of \gls{ADMM}.
The SCS solver \cite{O'Donoghue:2016} sets $\eps$ in Algorithm~\ref{alg:pcg} as a decreasing function of the \gls{ADMM} iteration counter $k$.
We adopt a different strategy in which $\eps$ is determined based on the \gls{ADMM} residuals.
In particular, we use
\[
  \eps \gets \max\left( \lambda\sqrt{\norm{\bar{r}^k_{\rm prim}}_\infty \norm{\bar{r}^k_{\rm dual}}_\infty}, \eps_{\rm min} \right),
\]
where $\bar{r}^k_{\rm prim}\eqdef E r^k_{\rm prim}$ and $\bar{r}^k_{\rm dual}\eqdef cD r^k_{\rm dual}$ are the scaled primal and dual residuals.
Parameter $\lambda\in ]0,1[$ ensures that $\eps$ is always lower than the geometric mean of the scaled primal and dual residuals.
We set $\lambda=0.15$ and $\eps_{\rm min}=10^{-7}$.
Since we use the $\ell_\infty$ norms for \gls{ADMM} residuals, we use the same norm in step~\ref{alg:pcg_termination} of Algorithm~\ref{alg:pcg}.
As $\eps$ depends on the \gls{ADMM} residuals, which are computed when evaluating \gls{ADMM} termination criteria, we evaluate these criteria after every $5$ \gls{ADMM} iterations.

\subsection{Matrix Equilibration}\label{subsec:matrix_equilibration}

\subsubsection{Computing column norms}
The \gls{CSR} representation of a sparse matrix allows for efficient computation of its row norms since the \texttt{RowPointer} array defines segments of the \texttt{Value} array corresponding to different rows.
Since we also store the matrix transpose, we can efficiently compute column norms of a matrix since they are equivalent to the row norms of its transpose.

A naive approach would be to have one thread per row computing its norm, but this approach is not the most efficient.
First, the workload may be distributed poorly among the threads since one row can have zero elements, and another can have many.
Second, the memory is accessed almost randomly as each thread iterates through its row, which can considerably deteriorate performance.

A more efficient way of computing the row norms of a matrix in the \gls{CSR} format is to represent the operation as a segmented reduction, where the segments are defined by the \texttt{RowPointer} array and, in the case of the $\ell_\infty$ norm, the associated binary operator is given by
\[
  x_1 \oplus x_2 = \max(|x_1|,|x_2|).
\]

\subsubsection{Matrix post-multiplication}
Matrix post-multiplication refers to evaluating the product $MD$, where $D\in\Re^{n\times n}$ is a diagonal matrix stored as an $n$-dimensional vector, and $M\in\Re^{m\times n}$ is a general sparse matrix in the \gls{CSR} format.
The \texttt{ColumnIndex} array can be used to determine the diagonal element of $D$ that multiplies each element of $A$,
\begin{lstlisting}[language=C]
  int column = ColumnIndex[idx];
  Value[idx] *= D[column];
\end{lstlisting}
This operation can be performed by many threads concurrently and independently.
As the memory read and write access to the array \texttt{Value} is fully coalesced, all memory addresses can be combined into a larger transaction.
However, the read access from \texttt{D} can be partly coalesced, but this does not impact the performance too much.

\subsubsection{Matrix pre-multiplication}
Matrix pre-multiplication in the product $DM$ is conceptually easier to implement since all elements in a row are multiplied with the same diagonal element of $D$.
However, it is not obvious how to determine the row index corresponding to an element of the \texttt{Value} array since the matrix $M$ is represented in the \gls{CSR} format.
We address this issue by computing \texttt{RowIndex} from the \texttt{RowPointer} array in advance, although it increases the memory usage.
The code that evaluates the matrix pre-multiplication is thus
\begin{lstlisting}[language=C]
  int row = RowIndex[idx];
  Value[idx] *= D[row];
\end{lstlisting}

\subsection{cuOSQP}

Table~\ref{tab:cuosqp} summarizes the main differences between OSQP and our GPU implementation of Algorithm~\ref{alg:osqp}.
Although our implementation requires two times more memory to store the problem matrices, it does not need to store any matrix factorizations.
Moreover, we can reduce the memory requirements by using the single-precision floating-point representation, which also leads to faster computations (see Section~\ref{sec:floats}).

We refer to our CUDA C implementation of the OSQP solver as cuOSQP.
The code is available online at
\begin{center}
  \url{https://github.com/oxfordcontrol/osqp/tree/cuda-1.0}
\end{center}
and its Python interface at
\begin{center}
  \url{https://github.com/oxfordcontrol/cuosqp}
\end{center}
cuOSQP uses cuBLAS, cuSPARSE, and Thrust libraries, which are included within the CUDA Toolkit.
Note that a custom implementation of linear algebra could improve efficiency of the solver even further.
However, relying on CUDA libraries not only saves the development time, but also ensures that our code is portable to various \glspl{GPU} and operating systems.
We have tested our code on both Linux and Windows machines, and have run it on a GeForce RTX 2080 Ti (launched in 2018) and a GeForce GTX 970 (launched in 2014).

\begin{table}[t]
	\centering
	\caption{The main differences between OSQP and cuOSQP implementations of \Alg~\ref{alg:osqp}.}
	\setlength{\tabcolsep}{1em}
	\begin{tabular}{lll}
    \toprule
    & OSQP (CPU) & cuOSQP (GPU) \\
    \midrule
		evaluating step~\ref{alg:solve_lin_sys} of \Alg~\ref{alg:osqp} & \tabitem linear system~\eqref{eqn:kkt} & \tabitem linear system~\eqref{eqn:kkt_red} \\
		& \tabitem $LDL^T$ factorization & \tabitem \gls{PCG} method (\Alg~\ref{alg:pcg}) \\
		\midrule
		parameter matrix $R$ & \tabitem $R = \diag(\rho_1,\ldots,\rho_m)$ & \tabitem $R=\bar{\rho}I$ \\
		& \tabitem $\rho_i$ set according to~\eqref{eqn:rho_vec} \\
		\midrule
		updating $\bar{\rho}$ & \tabitem rarely & \tabitem every $10$ iterations \\
		\midrule
		storing data matrices & \tabitem \gls{CSC} format & \tabitem \gls{CSR} format \\
		& \tabitem upper triangular $P$ & \tabitem full $P$ \\
		& \tabitem only $A$ & \tabitem both $A$ and $A\tpose$ \\
		\midrule
		checking termination & \tabitem every $25$ iterations & \tabitem every $5$ iterations \\
    \bottomrule
	\end{tabular}
	\label{tab:cuosqp}
\end{table}

\section{Numerical Results}\label{sec:numerics}

We evaluate performance of cuOSQP and compare it against both single- and multi-threaded versions of OSQP (version 0.6.0), which was shown to be competitive to and even faster than commercial \gls{QP} solvers \cite{Stellato:2020}.
Our main goal is to demonstrate how a parallel \gls{GPU} implementation can improve performance of an optimization solver for large-scale problems.
The sizes of benchmark problems range from $10^4$ to $10^8$ nonzero elements in $P$ and $A$.
We use the default parameters for both solvers.
By default, we use the single-precision floating-point representation with cuOSQP, but we also compare the single- and double-precision variants in Section~\ref{sec:floats}.

All numerical tests were performed on a Linux-based system with an i9-9900K @ 3.6GHz (8 cores) processor and 64 GB of DDR4 3200Mhz RAM, which is equipped with the NVIDIA GeForce RTX 2080 Ti \gls{GPU} with 11 GB of VRAM.

\subsection{OSQP Benchmark Problems}
We use the set of benchmark problems described in \cite[\App~A]{Stellato:2020}, which consist of \glspl{QP} from $7$ problem classes, ranging from standard random problems to applications in control, finance, statistics, and machine learning.
The problems are available online at \cite{osqp_benchmarks} and are summarized in the sequel.

\begin{itemize}
  \item \emph{Control}.
  The problem of controlling a linear time-invariant dynamical system can be formulated as the following constrained finite-time optimal control problem:
  \[
    \MinProblem{}{\displaystyle x_T\tpose Q_T x_T + \sum_{t=0}^{T-1} x_t\tpose Q x_t + u_t\tpose R u_t}{x_0 = x_{\rm init} \\ & x_{t+1} = Ax_t + Bu_t \\ & -\underbar{$x$} \le x_t \le \overline{x} \\ & -\underbar{$u$} \le u_t \le \overline{u}.}
  \]

  \item \emph{Equality}.
  This class consists of the following equality-constrained \glspl{QP}:
  \[
    \MinProblem{}{\half x\tpose P x + q\tpose x}{Ax = b.}
  \]

  \item \emph{Huber}.
  \emph{Huber fitting} or the \emph{robust least-squares problem} performs linear regression under the assumption that there are outliers in the data.
  The problem can be written as
  \[
    \text{minimize} \quad \sum_{i=1}^m \phi_{\rm hub}(a_i\tpose x - b_i),
  \]
  where the Huber penalty function $\phi_{\rm hub}\colon\Re\to\Re$ is defined as
  \[
    \phi_{\rm hub}(u) \eqdef \begin{cases} u^2 & |u| \le M \\ M(2|u|-M) & \text{otherwise}. \end{cases}
  \]

  \item \emph{Lasso}.
  The \emph{least absolute shrinkage and selection operator (lasso)} is a well-known technique aiming to obtain a sparse solution to a linear regression problem by adding an $\ell_1$ regularization term in the objective.
  The problem can be formulated as
  \[
    \text{minimize} \quad \norm{Ax-b}_2^2 + \lambda \norm{x}_1.
  \]

  \item \emph{Portfolio}.
  \emph{Portfolio optimization} is a problem arising in finance that seeks to allocate assets in a way that maximizes the risk-adjusted return.
  The problem has the following form:
  \[
    \MaxProblem{}{\mu\tpose x - \gamma x\tpose \Sigma x}{\boldsymbol{1}\tpose x = 1 \\ & x\ge 0,}
  \]
  where $x\in\Re^n$ represents the portfolio, $\mu\in\Re^n$ the vector of expected returns, $\gamma>0$ the risk-aversion parameter, and $\Sigma\in\symm_+^n$ the risk covariance matrix.

  \item \emph{Random}.
  This class consists of the following \gls{QP} with randomly generated data:
  \[
    \MinProblem{}{\half x\tpose P x + q\tpose x}{l \le Ax \le u.}
  \]

  \item \emph{SVM}.
  \emph{\Gls{SVM}} problem seeks an affine function that approximately classifies two sets of points.
  The problem can be stated as
  \[
    \text{minimize} \quad x\tpose x + \lambda \sum_{i=1}^m \max(0, b_i a_i\tpose x + 1),
  \]
  where $b_i\in\{-1,+1\}$ is the set label and $a_i$ the vector of features for the $i$-th point.
\end{itemize}

All instances were obtained from realistic non-trivial random data.
For each problem class we generate $10$ different instances for $15$ dimensions giving a total of $1050$ problems.
As a performance metric, we use the average runtime across $10$ different problem instances of the same size.

Figures~\ref{fig:results-1}--\ref{fig:results-2} show the computation runtimes achieved by OSQP and cuOSQP.
The figures show that OSQP is faster than cuOSQP for problem sizes of the order up to $10^5$.
However, for larger problem instances cuOSQP is significantly faster.
Furthermore, the slope of the runtimes achieved by OSQP is approximately constant, whereas for cuOSQP it is flatter for smaller problems and increases for larger.
This behavior is expected since smaller problems cannot fully utilize the \gls{GPU}, and the kernel launch and data transfer latencies cannot be amortized.
Moreover, the main focus of cuOSQP is on large-scale problems and thus we have not optimized it for small problem sizes.

Figure~\ref{fig:iter} shows the number of \gls{ADMM} iterations needed to satisfy the termination condition~\eqref{eqn:termination_condition}.
Updating $\bar{\rho}$ every $10$ iterations helps decrease the total number of \gls{ADMM} iterations for the problem classes Equality, Lasso, and SVM, which explains the obtained speedup shown in Figures~\ref{fig:results-1}--\ref{fig:results-2}.
For the Control, Portfolio, and Random classes the benefit is not apparent, while for the Huber class updating $\bar{\rho}$ less frequently seems to work better; in fact, our numerical tests indicate that the smallest number of iterations is achieved when $\bar{\rho}$ is kept constant.

\subsection{QDLDL}

When compared to OSQP's default single-threaded linear system solver QDLDL \cite{qdldl}, the maximum speedups achieved by cuOSQP range from $15$ to $270$ times.
The largest reduction in runtime is achieved for the Equality class, where OSQP takes $6.5~\rm{min}$ to solve the largest problem instance, while cuOSQP solves it in $1.4~\rm{s}$.
The second largest reduction is achieved for the SVM class with a reduction from $5.6~\rm{min}$ to $3.2~\rm{s}$.

For some problems, one can observe that OSQP runtimes do not necessarily increase with the problem size.
This behavior comes from computing the permutation of the KKT matrix prior to its factorization, which is performed by the AMD routine \cite{Amestoy:2004}, whose runtimes do not depend only on the number of nonzero elements in the KKT matrix.

\subsection{MKL Pardiso}

Apart from its single-threaded QDLDL linear system solver, OSQP can be interfaced with Intel MKL Pardiso \cite{pardiso}, a multi-threaded parallel direct sparse solver.
By default, Intel MKL Pardiso uses the maximum number of CPU cores available, which results in its best performance~\cite{pardiso}; hence, in our numerical tests the solver uses $8$ cores.
Figures~\ref{fig:results-1}--\ref{fig:results-2} show that the computation runtimes increase monotonically with the problem size when using MKL Pardiso.
Also, for smaller problem sizes OSQP is faster when using QDLDL, but for larger problems using MKL Pardiso reduces its runtimes significantly.
However, the maximum ratio of runtimes achieved by OSQP and cuOSQP is still between $3.7$ and $57$ times, depending on the problem class.

\begin{figure}[t]
  \centering
  \footnotesize
  \begin{tikzpicture}
    \begin{axis}[%
      hide axis,
      xmin=10,	
      xmax=50,	
      ymin=0,		
      ymax=0.4,	
      legend cell align=left,
      legend columns=3,
      legend style={%
        nodes={scale=1},
        fill=white,
        fill opacity=1,
        draw opacity=1,
        text opacity=1,
        /tikz/every even column/.append style={column sep=3em}
      },
      every axis plot/.append style={very thick}
      ]
      \addlegendimage{blue, mark=*, mark options={scale=1, fill=blue}}
      \addlegendentry{OSQP (QDLDL)};
      \addlegendimage{magenta, mark=pentagon*, mark options={scale=1, fill=magenta}}
      \addlegendentry{OSQP (MKL Pardiso)};
      \addlegendimage{red, mark=square*, mark options={scale=1, fill=red}}
      \addlegendentry{cuOSQP};
      \addlegendimage{black, mark=diamond*, mark options={scale=1, fill=black}}
      \addlegendentry{Speedup over QDLDL};
      \addlegendimage{gray, mark=triangle*, mark options={scale=1, fill=gray}}
      \addlegendentry{Speedup over MKL Pardiso};
    \end{axis}
  \end{tikzpicture}
  \\[.5em]
  \begin{tabular}{rl}
    \begin{tikzpicture}
      \begin{loglogaxis}[%
        title={Control},
        title style={at={(.1,.8)},anchor=west,fill=white,text=black,draw=black},
        width=\mywidth,
        height=\myheight,
        xmin = 5e3,
        xmax = 2e8,
        grid=both,
        grid style={line width=.1pt, draw=gray!10},
        major grid style={line width=.2pt,draw=gray!50},
        xticklabels={,,},
        ]
        \addplot [blue, mark=*, mark options={scale=1, fill=blue}] table[x=N, y=OSQP_runtime, col sep=comma] {data/control.csv};
        \addplot [magenta, mark=pentagon*, mark options={scale=1, fill=magenta}] table [x=N, y=OSQP_Pardiso_runtime, col sep=comma] {data/control.csv};
        \addplot [red, mark=square*, mark options={scale=1, fill=red}] table[x=N, y=cuOSQP3_runtime, col sep=comma] {data/control.csv};
      \end{loglogaxis}
    \end{tikzpicture}
    &
    \begin{tikzpicture}
      \begin{semilogxaxis}[%
        width=\mywidth,
        height=\myheight,
        xmin = 5e3,
        xmax = 2e8,
        grid=both,
        grid style={line width=.1pt, draw=gray!10},
        major grid style={line width=.2pt,draw=gray!50},
        ylabel near ticks,
        yticklabel pos=right,
        xticklabels={,,},
        ]
        \addplot [black, mark=diamond*, mark options={scale=1, fill=black}] table [x=N, y=speedup3, col sep=comma] {data/control.csv};
        \addplot [gray, mark=triangle*, mark options={scale=1, fill=gray}] table [x=N, y=speedup3_Pardiso, col sep=comma] {data/control.csv};
      \end{semilogxaxis}
    \end{tikzpicture}
    \\[-.75em]
    \begin{tikzpicture}
      \begin{loglogaxis}[%
        title={Equality},
        title style={at={(.1,.8)},anchor=west,fill=white,text=black,draw=black},
        width=\mywidth,
        height=\myheight,
        xmin = 5e3,
        xmax = 2e8,
        ymax = 3e3,
        grid=both,
        grid style={line width=.1pt, draw=gray!10},
        major grid style={line width=.2pt,draw=gray!50},
        xticklabels={,,},
        ]
        \addplot [blue, mark=*, mark options={scale=1, fill=blue}] table[x=N, y=OSQP_runtime, col sep=comma] {data/equality.csv};
        \addplot [magenta, mark=pentagon*, mark options={scale=1, fill=magenta}] table [x=N, y=OSQP_Pardiso_runtime, col sep=comma] {data/equality.csv};
        \addplot [red, mark=square*, mark options={scale=1, fill=red}] table[x=N, y=cuOSQP3_runtime, col sep=comma] {data/equality.csv};
      \end{loglogaxis}
    \end{tikzpicture}
    &
    \begin{tikzpicture}
      \begin{semilogxaxis}[%
        width=\mywidth,
        height=\myheight,
        xmin = 5e3,
        xmax = 2e8,
        grid=both,
        grid style={line width=.1pt, draw=gray!10},
        major grid style={line width=.2pt,draw=gray!50},
        ylabel near ticks,
        yticklabel pos=right,
        xticklabels={,,},
        ]
        \addplot [black, mark=diamond*, mark options={scale=1, fill=black}] table [x=N, y=speedup3, col sep=comma] {data/equality.csv};
        \addplot [gray, mark=triangle*, mark options={scale=1, fill=gray}] table [x=N, y=speedup3_Pardiso, col sep=comma] {data/equality.csv};
        \draw (4e7,5) rectangle (1.5e8,25);
      \end{semilogxaxis}
      \begin{semilogxaxis}[%
        width=0.3\textwidth,
        height=0.12\textheight,
        xmin = 4e7,
        xmax = 1.5e8,
        ymin = 5,
        ymax = 25,
        grid=both,
        grid style={line width=.1pt, draw=gray!10},
        major grid style={line width=.2pt,draw=gray!50},
        yticklabel pos=right,
        ytick={5,15,25},
        xticklabels={,,},
        shift={(0.05\textwidth,0.1\textheight)},
        axis background/.style={fill=white},
        ]
        \addplot [gray, very thick, mark=triangle*, mark options={scale=3, fill=gray}] table [x=N, y=speedup3_Pardiso, col sep=comma] {data/equality.csv};
      \end{semilogxaxis}
    \end{tikzpicture}
    \\[-.75em]
    \begin{tikzpicture}
      \begin{loglogaxis}[%
        title={Huber},
        title style={at={(.1,.8)},anchor=west,fill=white,text=black,draw=black},
        width=\mywidth,
        height=\myheight,
        xmin = 5e3,
        xmax = 2e8,
        ymax = 2e2,
        grid=both,
        grid style={line width=.1pt, draw=gray!10},
        major grid style={line width=.2pt,draw=gray!50},
        xticklabels={,,},
        ]
        \addplot [blue, mark=*, mark options={scale=1, fill=blue}] table[x=N, y=OSQP_runtime, col sep=comma] {data/huber.csv};
        \addplot [magenta, mark=pentagon*, mark options={scale=1, fill=magenta}] table [x=N, y=OSQP_Pardiso_runtime, col sep=comma] {data/huber.csv};
        \addplot [red, mark=square*, mark options={scale=1, fill=red}] table[x=N, y=cuOSQP3_runtime, col sep=comma] {data/huber.csv};
      \end{loglogaxis}
    \end{tikzpicture}
    &
    \begin{tikzpicture}
      \begin{semilogxaxis}[%
        width=\mywidth,
        height=\myheight,
        xmin = 5e3,
        xmax = 2e8,
        grid=both,
        grid style={line width=.1pt, draw=gray!10},
        major grid style={line width=.2pt,draw=gray!50},
        ylabel near ticks,
        yticklabel pos=right,
        xticklabels={,,},
        ]
        \addplot [black, mark=diamond*, mark options={scale=1, fill=black}] table [x=N, y=speedup3, col sep=comma] {data/huber.csv};
        \addplot [gray, mark=triangle*, mark options={scale=1, fill=gray}] table [x=N, y=speedup3_Pardiso, col sep=comma] {data/huber.csv};
      \end{semilogxaxis}
    \end{tikzpicture}
    \\[-.75em]
    \begin{tikzpicture}
      \begin{loglogaxis}[%
        title={Lasso},
        title style={at={(.1,.8)},anchor=west,fill=white,text=black,draw=black},
        width=\mywidth,
        height=\myheight,
        xmin = 5e3,
        xmax = 2e8,
        grid=both,
        grid style={line width=.1pt, draw=gray!10},
        major grid style={line width=.2pt,draw=gray!50},
        legend pos = north west,
        legend cell align={left},
        xlabel={Problem size $N$},
        ]
        \addplot [blue, mark=*, mark options={scale=1, fill=blue}] table[x=N, y=OSQP_runtime, col sep=comma] {data/lasso.csv};
        \addplot [magenta, mark=pentagon*, mark options={scale=1, fill=magenta}] table [x=N, y=OSQP_Pardiso_runtime, col sep=comma] {data/lasso.csv};
        \addplot [red, mark=square*, mark options={scale=1, fill=red}] table[x=N, y=cuOSQP3_runtime, col sep=comma] {data/lasso.csv};
      \end{loglogaxis}
    \end{tikzpicture}
    &
    \begin{tikzpicture}
      \begin{semilogxaxis}[%
        width=\mywidth,
        height=\myheight,
        xmin = 5e3,
        xmax = 2e8,
        grid=both,
        grid style={line width=.1pt, draw=gray!10},
        major grid style={line width=.2pt,draw=gray!50},
        legend pos = north west,
        legend cell align={left},
        ylabel near ticks,
        yticklabel pos=right,
        xlabel={Problem size $N$},
        ]
        \addplot [black, mark=diamond*, mark options={scale=1, fill=black}] table [x=N, y=speedup3, col sep=comma] {data/lasso.csv};
        \addplot [gray, mark=triangle*, mark options={scale=1, fill=gray}] table [x=N, y=speedup3_Pardiso, col sep=comma] {data/lasso.csv};
      \end{semilogxaxis}
    \end{tikzpicture}
  \end{tabular}
  \caption{%
    Numerical comparison of OSQP and cuOSQP for problem classes Control, Equality, Huber, and Lasso.
    Left: Average computation runtimes (in seconds) across $10$ problem instances as a function of problem size~$N$.
    Right: Ratio between the runtimes achieved with OSQP and cuOSQP.
  }
  \label{fig:results-1}
\end{figure}

\begin{figure}[t]
  \centering
  \footnotesize
  \begin{tikzpicture}
    \begin{axis}[%
      hide axis,
      xmin=10,	
      xmax=50,	
      ymin=0,		
      ymax=0.4,	
      legend cell align=left,
      legend columns=3,
      legend style={%
        nodes={scale=1},
        fill=white,
        fill opacity=1,
        draw opacity=1,
        text opacity=1,
        /tikz/every even column/.append style={column sep=3em}
      },
      every axis plot/.append style={very thick}
      ]
      \addlegendimage{blue, mark=*, mark options={scale=1, fill=blue}}
      \addlegendentry{OSQP (QDLDL)};
      \addlegendimage{magenta, mark=pentagon*, mark options={scale=1, fill=magenta}}
      \addlegendentry{OSQP (MKL Pardiso)};
      \addlegendimage{red, mark=square*, mark options={scale=1, fill=red}}
      \addlegendentry{cuOSQP};
      \addlegendimage{black, mark=diamond*, mark options={scale=1, fill=black}}
      \addlegendentry{Speedup over QDLDL};
      \addlegendimage{gray, mark=triangle*, mark options={scale=1, fill=gray}}
      \addlegendentry{Speedup over MKL Pardiso};
    \end{axis}
  \end{tikzpicture}
  \\[.5em]
  \begin{tabular}{rl}
    \begin{tikzpicture}
      \begin{loglogaxis}[%
        title={Portfolio},
        title style={at={(.1,.8)},anchor=west,fill=white,text=black,draw=black},
        width=\mywidth,
        height=\myheight,
        xmin = 5e3,
        xmax = 2e8,
        grid=both,
        grid style={line width=.1pt, draw=gray!10},
        major grid style={line width=.2pt,draw=gray!50},
        xticklabels={,,},
        ]
        \addplot [blue, mark=*, mark options={scale=1, fill=blue}] table[x=N, y=OSQP_runtime, col sep=comma] {data/portfolio.csv};
        \addplot [magenta, mark=pentagon*, mark options={scale=1, fill=magenta}] table [x=N, y=OSQP_Pardiso_runtime, col sep=comma] {data/portfolio.csv};
        \addplot [red, mark=square*, mark options={scale=1, fill=red}] table[x=N, y=cuOSQP3_runtime, col sep=comma] {data/portfolio.csv};
      \end{loglogaxis}
    \end{tikzpicture}
    &
    \begin{tikzpicture}
      \begin{semilogxaxis}[%
        width=\mywidth,
        height=\myheight,
        xmin = 5e3,
        xmax = 2e8,
        grid=both,
        grid style={line width=.1pt, draw=gray!10},
        major grid style={line width=.2pt,draw=gray!50},
        ylabel near ticks,
        yticklabel pos=right,
        xticklabels={,,},
        ]
        \addplot [black, mark=diamond*, mark options={scale=1, fill=black}] table [x=N, y=speedup3, col sep=comma] {data/portfolio.csv};
        \addplot [gray, mark=triangle*, mark options={scale=1, fill=gray}] table [x=N, y=speedup3_Pardiso, col sep=comma] {data/portfolio.csv};
      \end{semilogxaxis}
    \end{tikzpicture}
    \\[-.75em]
    \begin{tikzpicture}
      \begin{loglogaxis}[%
        title={Random},
        title style={at={(.1,.8)},anchor=west,fill=white,text=black,draw=black},
        width=\mywidth,
        height=\myheight,
        xmin = 5e3,
        xmax = 2e8,
        grid=both,
        grid style={line width=.1pt, draw=gray!10},
        major grid style={line width=.2pt,draw=gray!50},
        xticklabels={,,},
        ]
        \addplot [blue, mark=*, mark options={scale=1, fill=blue}] table[x=N, y=OSQP_runtime, col sep=comma] {data/random.csv};
        \addplot [magenta, mark=pentagon*, mark options={scale=1, fill=magenta}] table [x=N, y=OSQP_Pardiso_runtime, col sep=comma] {data/random.csv};
        \addplot [red, mark=square*, mark options={scale=1, fill=red}] table[x=N, y=cuOSQP3_runtime, col sep=comma] {data/random.csv};
      \end{loglogaxis}
    \end{tikzpicture}
    &
    \begin{tikzpicture}
      \begin{semilogxaxis}[%
        width=\mywidth,
        height=\myheight,
        xmin = 5e3,
        xmax = 2e8,
        grid=both,
        grid style={line width=.1pt, draw=gray!10},
        major grid style={line width=.2pt,draw=gray!50},
        ylabel near ticks,
        yticklabel pos=right,
        xticklabels={,,},
        ]
        \addplot [black, mark=diamond*, mark options={scale=1, fill=black}] table [x=N, y=speedup3, col sep=comma] {data/random.csv};
        \addplot [gray, mark=triangle*, mark options={scale=1, fill=gray}] table [x=N, y=speedup3_Pardiso, col sep=comma] {data/random.csv};
      \end{semilogxaxis}
    \end{tikzpicture}
    \\[-.75em]
    \begin{tikzpicture}
      \begin{loglogaxis}[%
        title={SVM},
        title style={at={(.1,.8)},anchor=west,fill=white,text=black,draw=black},
        width=\mywidth,
        height=\myheight,
        xmin = 5e3,
        xmax = 2e8,
        grid=both,
        grid style={line width=.1pt, draw=gray!10},
        major grid style={line width=.2pt,draw=gray!50},
        legend pos = north west,
        legend cell align={left},
        xlabel={Problem size $N$},
        ]
        \addplot [blue, mark=*, mark options={scale=1, fill=blue}] table[x=N, y=OSQP_runtime, col sep=comma] {data/svm.csv};
        \addplot [magenta, mark=pentagon*, mark options={scale=1, fill=magenta}] table [x=N, y=OSQP_Pardiso_runtime, col sep=comma] {data/svm.csv};
        \addplot [red, mark=square*, mark options={scale=1, fill=red}] table[x=N, y=cuOSQP3_runtime, col sep=comma] {data/svm.csv};
      \end{loglogaxis}
    \end{tikzpicture}
    &
    \begin{tikzpicture}
      \begin{semilogxaxis}[%
        width=\mywidth,
        height=\myheight,
        xmin = 5e3,
        xmax = 2e8,
        grid=both,
        grid style={line width=.1pt, draw=gray!10},
        major grid style={line width=.2pt,draw=gray!50},
        ylabel near ticks,
        yticklabel pos=right,
        xlabel={Problem size $N$},
        ]
        \addplot [black, mark=diamond*, mark options={scale=1, fill=black}] table [x=N, y=speedup3, col sep=comma] {data/svm.csv};
        \addplot [gray, mark=triangle*, mark options={scale=1, fill=gray}] table [x=N, y=speedup3_Pardiso, col sep=comma] {data/svm.csv};
      \end{semilogxaxis}
    \end{tikzpicture}
  \end{tabular}
  \caption{%
    Numerical comparison of OSQP and cuOSQP for problem classes Portfolio, Random, and SVM.
    Left: Average computation runtimes (in seconds) across $10$ problem instances as a function of problem size~$N$.
    Right: Ratio between the runtimes achieved with OSQP and cuOSQP.
  }
  \label{fig:results-2}
\end{figure}

\begin{figure}[t]
  \centering
  \footnotesize
  \begin{tikzpicture}
    \begin{axis}[%
      hide axis,
      xmin=10,	
      xmax=50,	
      ymin=0,		
      ymax=0.4,	
      legend cell align=left,
      legend columns=2,
      legend style={%
        nodes={scale=1},
        fill=white,
        fill opacity=1,
        draw opacity=1,
        text opacity=1,
        /tikz/every even column/.append style={column sep=3em}
      },
      every axis plot/.append style={very thick}
      ]
      \addlegendimage{blue, mark=*, mark options={scale=1, fill=blue}}
      \addlegendentry{OSQP};
      \addlegendimage{red, mark=square*, mark options={scale=1, fill=red}}
      \addlegendentry{cuOSQP};
    \end{axis}
  \end{tikzpicture}
  \\[.5em]
  \begin{tabular}[t]{rl}
    \begin{tikzpicture}
      \begin{semilogxaxis}[%
        title={Control},
        title style={at={(.1,.8)},anchor=west,fill=white,text=black,draw=black},
        width=\mywidth,
        height=\myheight,
        xmin = 5e3,
        xmax = 2e8,
        grid=both,
        grid style={line width=.1pt, draw=gray!10},
        major grid style={line width=.2pt,draw=gray!50},
        xticklabels={,,},
        ]
        \addplot [blue, mark=*, mark options={scale=1, fill=blue}, error bars/.cd, y explicit, y dir=both] table[x=N, y=OSQP_iter_mean, y error=OSQP_iter_std, col sep=comma] {data/control.csv};
        \addplot[red, mark=square*, mark options={scale=1, fill=red}, error bars/.cd, y explicit, y dir=both] table [x=N, y=cuOSQP3_iter_mean, y error=cuOSQP3_iter_std, col sep=comma] {data/control.csv};
      \end{semilogxaxis}
    \end{tikzpicture}
    &
    \begin{tikzpicture}
      \begin{semilogxaxis}[%
        title={Equality},
        title style={at={(.1,.8)},anchor=west,fill=white,text=black,draw=black},
        width=\mywidth,
        height=\myheight,
        xmin = 5e3,
        xmax = 2e8,
        grid=both,
        grid style={line width=.1pt, draw=gray!10},
        major grid style={line width=.2pt,draw=gray!50},
        xticklabels={,,},
        yticklabel pos=right,
        ]
        \addplot [blue, mark=*, mark options={scale=1, fill=blue}, error bars/.cd, y explicit, y dir=both] table[x=N, y=OSQP_iter_mean, y error=OSQP_iter_std, col sep=comma] {data/equality.csv};
        \addplot[red, mark=square*, mark options={scale=1, fill=red}, error bars/.cd, y explicit, y dir=both] table [x=N, y=cuOSQP3_iter_mean, y error=cuOSQP3_iter_std, col sep=comma] {data/equality.csv};
      \end{semilogxaxis}
    \end{tikzpicture}
    \\[-.75em]
    \begin{tikzpicture}
      \begin{semilogxaxis}[%
        title={Huber},
        title style={at={(.1,.8)},anchor=west,fill=white,text=black,draw=black},
        width=\mywidth,
        height=\myheight,
        xmin = 5e3,
        xmax = 2e8,
        grid=both,
        grid style={line width=.1pt, draw=gray!10},
        major grid style={line width=.2pt,draw=gray!50},
        xticklabels={,,},
        ]
        \addplot [blue, mark=*, mark options={scale=1, fill=blue}, error bars/.cd, y explicit, y dir=both] table[x=N, y=OSQP_iter_mean, y error=OSQP_iter_std, col sep=comma] {data/huber.csv};
        \addplot[red, mark=square*, mark options={scale=1, fill=red}, error bars/.cd, y explicit, y dir=both] table [x=N, y=cuOSQP3_iter_mean, y error=cuOSQP3_iter_std, col sep=comma] {data/huber.csv};
      \end{semilogxaxis}
    \end{tikzpicture}
    &
    \begin{tikzpicture}
      \begin{semilogxaxis}[%
        title={Lasso},
        title style={at={(.1,.8)},anchor=west,fill=white,text=black,draw=black},
        width=\mywidth,
        height=\myheight,
        xmin = 5e3,
        xmax = 2e8,
        grid=both,
        grid style={line width=.1pt, draw=gray!10},
        major grid style={line width=.2pt,draw=gray!50},
        xticklabels={,,},
        yticklabel pos=right,
        ]
        \addplot [blue, mark=*, mark options={scale=1, fill=blue}, error bars/.cd, y explicit, y dir=both] table[x=N, y=OSQP_iter_mean, y error=OSQP_iter_std, col sep=comma] {data/lasso.csv};
        \addplot[red, mark=square*, mark options={scale=1, fill=red}, error bars/.cd, y explicit, y dir=both] table [x=N, y=cuOSQP3_iter_mean, y error=cuOSQP3_iter_std, col sep=comma] {data/lasso.csv};
      \end{semilogxaxis}
    \end{tikzpicture}
    \\[-.75em]
    \begin{tikzpicture}
      \begin{semilogxaxis}[%
        title={Portfolio},
        title style={at={(.1,.8)},anchor=west,fill=white,text=black,draw=black},
        width=\mywidth,
        height=\myheight,
        xmin = 5e3,
        xmax = 2e8,
        ymax = 450,
        grid=both,
        grid style={line width=.1pt, draw=gray!10},
        major grid style={line width=.2pt,draw=gray!50},
        xticklabels={,,},
        xlabel={\phantom{Problem size $N$}},
        ]
        \addplot [blue, mark=*, mark options={scale=1, fill=blue}, error bars/.cd, y explicit, y dir=both] table[x=N, y=OSQP_iter_mean, y error=OSQP_iter_std, col sep=comma] {data/portfolio.csv};
        \addplot[red, mark=square*, mark options={scale=1, fill=red}, error bars/.cd, y explicit, y dir=both] table [x=N, y=cuOSQP3_iter_mean, y error=cuOSQP3_iter_std, col sep=comma] {data/portfolio.csv};
      \end{semilogxaxis}
    \end{tikzpicture}
    &
    \begin{tikzpicture}
      \begin{semilogxaxis}[%
        title={Random},
        title style={at={(.1,.8)},anchor=west,fill=white,text=black,draw=black},
        width=\mywidth,
        height=\myheight,
        xmin = 5e3,
        xmax = 2e8,
        ymax = 1.15e3,
        grid=both,
        grid style={line width=.1pt, draw=gray!10},
        major grid style={line width=.2pt,draw=gray!50},
        xticklabels={,,},
        yticklabel pos=right,
        xlabel={Problem size $N$},
        ]
        \addplot [blue, mark=*, mark options={scale=1, fill=blue}, error bars/.cd, y explicit, y dir=both] table[x=N, y=OSQP_iter_mean, y error=OSQP_iter_std, col sep=comma] {data/random.csv};
        \addplot[red, mark=square*, mark options={scale=1, fill=red}, error bars/.cd, y explicit, y dir=both] table [x=N, y=cuOSQP3_iter_mean, y error=cuOSQP3_iter_std, col sep=comma] {data/random.csv};
      \end{semilogxaxis}
    \end{tikzpicture}
    \\[-1.3em]
    \begin{tikzpicture}
      \begin{semilogxaxis}[%
        title={SVM},
        title style={at={(.1,.8)},anchor=west,fill=white,text=black,draw=black},
        width=\mywidth,
        height=\myheight,
        xmin = 5e3,
        xmax = 2e8,
        grid=both,
        grid style={line width=.1pt, draw=gray!10},
        major grid style={line width=.2pt,draw=gray!50},
        xticklabels={,,},
        xlabel={Problem size $N$},
        ]
        \addplot [blue, mark=*, mark options={scale=1, fill=blue}, error bars/.cd, y explicit, y dir=both] table[x=N, y=OSQP_iter_mean, y error=OSQP_iter_std, col sep=comma] {data/svm.csv};
        \addplot[red, mark=square*, mark options={scale=1, fill=red}, error bars/.cd, y explicit, y dir=both] table [x=N, y=cuOSQP3_iter_mean, y error=cuOSQP3_iter_std, col sep=comma] {data/svm.csv};
      \end{semilogxaxis}
    \end{tikzpicture}
  \end{tabular}
  \caption{%
    Number of ADMM iterations needed to reach a termination criterion by OSQP and cuOSQP for the $7$ problem classes.
    The markers show the average number of iterations across $10$ problem instances as a function of problem size~$N$, while the vertical bars show the standard deviation.
  }
  \label{fig:iter}
\end{figure}

\subsection{Floating-Point Precision}\label{sec:floats}

Figure~\ref{fig:floats} shows the average computation times when running cuOSQP on the Portfolio benchmark class for both single- and double-precision floating-point representations.
The penalty in computation times when using double- over single-precision is less than $2$ times over all problem sizes.
Moreover, our numerical results suggest that for other problem classes this penalty is even smaller (data not shown).
This is counter-intuitive at first since the \gls{GPU} used in our tests has $32$ times higher single-precision floating-point performance than in double-precision.
However, most numerical methods that we use, especially \gls{SpMV}, are memory-bound operations, which means that the computation times are limited by the memory bandwidth.
Hence, we expect that the achieved speedups would be even larger for \glspl{GPU} with higher memory bandwidths, such as NVIDIA V100 or V100s models.

\begin{figure}[t]
  \centering
  \footnotesize
  \begin{tikzpicture}
    \begin{axis}[%
      hide axis,
      xmin=10,	
      xmax=50,	
      ymin=0,		
      ymax=0.4,	
      legend cell align=left,
      legend columns=3,
      legend style={%
        nodes={scale=1},
        fill=white,
        fill opacity=1,
        draw opacity=1,
        text opacity=1,
        /tikz/every even column/.append style={column sep=3em}
      },
      every axis plot/.append style={very thick}
      ]
      \addlegendimage{red, mark=square*, mark options={scale=1, fill=red}}
      \addlegendentry{cuOSQP (\texttt{float})};
      \addlegendimage{brown, mark=triangle*, mark options={scale=1, fill=brown}}
      \addlegendentry{cuOSQP (\texttt{double})};
      \addlegendimage{black, mark=diamond*, mark options={scale=1, fill=black}}
      \addlegendentry{Speedup};
    \end{axis}
  \end{tikzpicture}
  \\[.5em]
  \begin{tabular}{rl}
    \begin{tikzpicture}
      \begin{loglogaxis}[%
        title={Portfolio},
        title style={at={(.1,.8)},anchor=west,fill=white,text=black,draw=black},
        width=\mywidth,
        height=\myheight,
        xlabel={Problem size $N$},
        grid=both,
        grid style={line width=.1pt, draw=gray!10},
        major grid style={line width=.2pt,draw=gray!50},
        ]
        \addplot [brown, mark=triangle*, mark options={scale=1, fill=brown}] table[x=N, y=cuOSQP3_double_runtime, col sep=comma] {data/portfolio.csv};
        \addplot [red, mark=square*, mark options={scale=1, fill=red}] table[x=N, y=cuOSQP3_runtime, col sep=comma] {data/portfolio.csv};
      \end{loglogaxis}
    \end{tikzpicture}
    &
    \begin{tikzpicture}
      \begin{semilogxaxis}[%
        width=\mywidth,
        height=\myheight,
        xlabel={Problem size $N$},
        grid=both,
        grid style={line width=.1pt, draw=gray!10},
        major grid style={line width=.2pt,draw=gray!50},
        ylabel near ticks,
        yticklabel pos=right,
        ]
        \addplot [black, mark=diamond*, mark options={scale=1, fill=black}] table [x=N, y=speedup3_double, col sep=comma] {data/portfolio.csv};
      \end{semilogxaxis}
    \end{tikzpicture}
  \end{tabular}
  \caption{%
    Numerical comparison of cuOSQP compiled with single- and double-precision floating-point representation for the Portfolio problem class.
    Left: Average computation runtimes (in seconds) across $10$ problem instances as a function of problem size $N$.
    Right: Ratio between the achieved runtimes.
  }
  \label{fig:floats}
\end{figure}

\section{Conclusions}\label{sec:conclusion}

We have explored the possibilities offered by the massive parallelism of \glspl{GPU} to accelerate solutions to large-scale \glspl{QP} and have managed to solve problems with hundreds of millions nonzero entries in the problem matrices in only a few seconds.
Our implementation cuOSQP is built on top of OSQP, a state-of-the-art \gls{QP} solver based on \gls{ADMM}.
The large speedup is achieved by using the \gls{PCG} method for solving the linear system arising in ADMM and by parallelizing all vector and matrix operations.
Our numerical tests confirm that GPUs are not suited for solving small problems for which the CPU implementation is generally much faster.
Our open-source implementation is written in CUDA C, and has been tested on both Linux and Windows machines.

Our implementation stores all problem data and ADMM iterates in the GPU memory.
While this design choice reduces the size of memory transferred between the system and the GPU, its drawback is that the size of problems are limited by the available GPU memory.
One possible extension would be to use the \emph{unified memory} approach, which merges the system memory with the GPU memory, and automatically transfers data on demand between the two memory spaces.
Alternatively, we could use multiple \glspl{GPU} to solve problems whose data could not fit on a single GPU.
The main challenges with a multi-GPU approach include the distribution of the workload across multiple devices and ensuring synchronization between them.

\section*{Acknowledgements}

We are grateful to Samuel Balula for helpful discussions and managing the hardware used in this work.
This project has received funding from the European Research Council (ERC) under the European Union's Horizon 2020 research and innovation programme grant agreement OCAL, No.\ 787845.

\bibliography{refs}

\end{document}

%% file: defs.tex
\usepackage{abstract} 
\usepackage{parskip}

\usepackage{titlesec}
\def\sectionfont{\sffamily\Large\bfseries\boldmath}
\def\subsectionfont{\sffamily\large\bfseries\boldmath}
\def\paragraphfont{\sffamily\normalsize\bfseries\boldmath}
\titleformat*{\section}{\sectionfont}
\titleformat*{\subsection}{\subsectionfont}
\titleformat*{\subsubsection}{\paragraphfont}
\titleformat*{\paragraph}{\paragraphfont}
\titleformat*{\subparagraph}{\paragraphfont}
\usepackage[small,labelfont={bf,sf}]{caption}  
\usepackage{fullpage,graphicx,psfrag,url}
\usepackage{multirow}
\usepackage{enumitem}
\setlist{nolistsep}

\usepackage{booktabs}   
\usepackage{csvsimple}	
\usepackage{adjustbox}  

\newcounter{algorithmctr}[section]
\renewcommand{\thealgorithmctr}{\thesection.\arabic{algorithmctr}}

\usepackage{amsmath,amsthm}

\newtheoremstyle{exampstyle}
  {.5\baselineskip} 
  {\topsep} 
  {} 
  {} 
  {\bfseries} 
  {.} 
  {.5em} 
  {} 
\theoremstyle{exampstyle}

\usepackage[english]{babel}
\usepackage{algorithm}
\usepackage[noend]{algpseudocode}
\usepackage{amssymb,mathtools,amsfonts}
\algrenewcommand\algorithmicindent{1em}%

\usepackage{hyperref}
\hypersetup{%
    colorlinks=true, 
    linktoc=all,     
    linkcolor=black,  
    citecolor= black,
    urlcolor=blue
}

\usepackage[xindy, acronyms]{glossaries}   
\glsdisablehyper

\newcommand{\ie}{\emph{i.e.}, }

\newcommand{\mcf}{\mathcal}

\renewcommand{\Re}{{\mbox{\bf R}}}

\newcommand{\symm}{{\mbox{\bf S}}}

\newcommand{\eqdef}{\coloneqq}
\newcommand{\eps}{\varepsilon}

\newcommand{\Alg}{{Alg.}}
\newcommand{\App}{{Appendix}}

\newcommand{\Prop}{{Prop.}}
\newcommand{\Thm}{{Thm.}}

\DeclareMathOperator*{\argmin}{\operatorname{argmin}}
\DeclareMathOperator{\diag}{diag}

\DeclareMathOperator{\mean}{mean}
\DeclareMathOperator{\nnz}{nnz}

\newcommand{\project}[1]{{\rm\Pi}_{#1}}

\newcommand{\tpose}{^T}
\newcommand{\inv}{^{-1}}

\newcommand{\norm}[1]{\lVert#1\rVert}

\newcommand{\half}{\tfrac{1}{2}}

\usepackage[binary-units=true]{siunitx}
\usepackage{tikz} 
\usepackage{pgfplots}
\pgfplotsset{compat=newest} 
\usepgfplotslibrary{units} 
\pgfplotsset{%
    /pgfplots/ybar legend/.style={
    /pgfplots/legend image code/.code={%
       \draw[##1,/tikz/.cd,
       yshift=-0.3em
       ]
        (0cm,0cm) rectangle (5pt,0.8em);},
   }
}
\pgfplotsset{select coords between index/.style 2 args={
    x filter/.code={
        \ifnum\coordindex<#1\fi
        \ifnum\coordindex>#2\fi
    }
}}

\newcommand{\Cpp}{{C\nolinebreak[4]\hspace{-.05em}\raisebox{.4ex}{\tiny\bf ++}}\ }

\newacronym{ADMM}{ADMM}{alternating direction method of multipliers}
\newacronym{BLAS}{BLAS}{Basic Linear Algebra Subprograms}
\newacronym{CG}{CG}{conjugate gradient}
\newacronym{COO}{COO}{coordinate}
\newacronym{CSR}{CSR}{compressed sparse row}
\newacronym{CSC}{CSC}{compressed sparse column}
\newacronym{CUDA}{CUDA}{Compute Unified Device Architecture}
\newacronym{GPU}{GPU}{graphics processing unit}
\newacronym{LP}{LP}{linear program}
\newacronym{PCG}{PCG}{preconditioned conjugate gradient}
\newacronym{QP}{QP}{quadratic program}
\newacronym{SIMD}{SIMD}{single instruction, multiple data}
\newacronym{SM}{SM}{streaming multiprocessor}
\newacronym{SpMV}{SpMV}{Sparse Matrix-Vector multiplication}
\newacronym{STL}{STL}{Standard Template Library}
\newacronym{SVM}{SVM}{support vector machine}

\newcommand{\OptProblem}[5][]{
\begin{aligned}\label{#1}
  \begin{array}{ll}
    \underset{#3}{\rm{#2}}&\hspace{-0ex}#4\vspace{.5ex}\\
    \rm{subject~to}&#5
  \end{array}
\end{aligned}
}

\newcommand{\MinProblem}[4][]
{\OptProblem[#1]{minimize}{#2}{#3}{#4}}
\newcommand{\MaxProblem}[4][]
{\OptProblem[#1]{maximize}{#2}{#3}{#4}}

\usepackage{xcolor}

\usepackage{listings}
\definecolor{codegreen}{rgb}{0,0.6,0}
\definecolor{codegray}{rgb}{0.5,0.5,0.5}
\definecolor{codepurple}{rgb}{0.58,0,0.82}
\definecolor{backcolor}{rgb}{0.95,0.95,0.92}
\lstdefinestyle{mystyle}{
    backgroundcolor=\color{backcolor},
    commentstyle=\color{codegreen},
    keywordstyle=\color{magenta},
    numberstyle=\tiny\color{codegray},
    stringstyle=\color{codepurple},
    basicstyle=\ttfamily\footnotesize,
    breakatwhitespace=false,
    breaklines=true,
    captionpos=b,
    keepspaces=true,
    showspaces=false,
    showstringspaces=false,
    showtabs=false,
    tabsize=2
}